\documentclass[leqno,12pt,twoside]{amsart}
\usepackage{amsmath,amscd,amssymb,upgreek} 
\usepackage[applemac]{inputenc}
\usepackage[all]{xy}
\usepackage[english]{babel}
\usepackage{graphicx}
\newtheorem{Lem}{Lemma \thesection.\!\!}
\newtheorem{Th}[Lem]{Theorem \thesection.\!\!}
\newtheorem{Cor}[Lem]{Corollary \thesection.\!\!}
\newtheorem{Def}[Lem]{Definition \thesection.\!\!}

\newtheorem{Prop}[Lem]{Proposition \thesection.\!\!}

\newtheorem{Lem and Def}[Lem]{Lemma and Definition \thesection.\!\!}

\newtheorem{Prop and Def}[Lem]{Proposition and Definition   \thesection.\!\!}

 \def\cal{\mathcal}
 
\def\a{\alpha }
\def\b{\beta }

\def\g{\gamma }

\def\l{\lambda }
\def\m{\mu }

\def\s{\sigma }
\def\t{\theta }

\def\a{\mathfrak a}  
\def\g{\mathfrak g}
 
\def\m{\mathfrak m}
\def\n{\mathfrak n }
\def\h{\mathfrak h } 
\def\b{\mathfrak b } 
\def\l{\mathfrak l } 
\def\s{\mathfrak s } 
\def\t{\mathfrak t } 
\def\z{\mathfrak z}
\def\tl{\triangleleft}
  
\def\cal{\mathcal}

\newcommand{\C}{\mathbb{C}}

\newcommand{\R}{\mathbb{R}}

\newcommand{\Z}{\mathbb{Z}}
\newcommand{\ra}{\rightarrow}

\begin{document}
\title[Homogeneous CR-Solvmanifolds]{Homogeneous CR-Solvmanifolds as K\"ahler Obstructions}
\author{Bruce GILLIGAN and Karl OELJEKLAUS} 
\date{\today}

\thanks{This work was partially done during a visit of the first author at the LATP at 
the Universit\'{e} de Provence in June 2008 and 
of the second author at the University of Regina in October-November 2008. 
We would like to thank the respective institutions for their hospitality and financial support. 
This work was also partially supported by an NSERC Discovery Grant.}

\address{{\it Bruce Gilligan}: 
Department of Mathematics and Statistics,   
University of Regina,    
3737 Wascana Boulevard, 
Regina, Canada S4S 0A2.  }
\email{gilligan@math.uregina.ca }
 
\address{{\it Karl Oeljeklaus}:
LATP-UMR(CNRS) 6632, 
CMI-Universit\'e d'Aix-Marseille I, 39, rue Joliot-Curie, F-13453 Marseille Cedex 13, France.}
\email{karloelj@cmi.univ-mrs.fr}

\begin{abstract}
We give a precise characterization when a compact homogeneous CR-solvmanifold
is CR-embeddable in a complex K\"ahler manifold. Equivalently this gives
a non-K\"ahler criterion for complex manifolds containing CR-solvmanifolds
not satisfying these conditions.
This paper is the natural continuation of \cite{OR} and \cite{GOR}.
\end{abstract}

\setcounter{section}{0}
\noindent
\maketitle

\section{Introduction}  

There are many results known about the structure of real solvmanifolds.   
One of these is the conjecture of Mostow \cite{Mos}, subsequently proved by L. Auslander, e.g., see \cite{Aus} for the history, 
that every solvmanifold is a vector bundle over a compact solvmanifold.  
In the category of complex solvmanifolds one would additionally like to understand the structure of 
the manifold with respect to complex analytic objects defined on it and, particularly, the role played by the base of 
the vector bundle noted above - this base determines the topology, so does it also control the complex analysis?   
Because of the connection with the existence of plurisubharmonic functions and analytic 
hypersurfaces, one problem of this type concerns the existence of a K\"{a}hler metric.  

\vskip 2ex\noindent  
 In the fundamental paper \cite{Lo}, triples $(\Gamma, G_0, G)$ are considered, 
where $G$ is a connected simply-connected solvable complex Lie group, 
$G_0$ is a connected (simply-connected) {\it totally real } subgroup of $G$, 
and $\Gamma$ is discrete and cocompact in $G_0$.  It is proven that
the complex manifold $X:= G/\Gamma$
admits a K\"ahler metric if and only if
for every element $x \in \g_0$ in the Lie algebra of $G_0$ the operator
${\rm ad} x$ has {\it purely imaginary spectrum} on $\g_0$.\\

\vskip 2ex\noindent 
In \cite {OR} the more general situation
of triples $(\Gamma, G_0, G)$ 
with Lie algebras $\g_0$ and $\g$ like above and
$\mathfrak{g}_0 + i \mathfrak{g}_0 =\mathfrak{g}$, but not necessarily
totally real $\g_0$, is studied.
Note that 
\[      
        \mathfrak{m}:=\mathfrak{g}_0 \cap i \mathfrak{g}_0    
\]  
is an ideal in $\mathfrak{g}$ and therefore also in $\mathfrak{g}_0$.   
Amongst other things, it is shown in \cite{OR} that $G/\Gamma$ being K\"{a}hler 
implies that the adjoint action of $\mathfrak{g}_0$ on $ \mathfrak{m}$ has purely imaginary spectrum.  

\vskip 2ex\noindent 
We underline that in both papers
\cite{Lo} and \cite{OR} the K\"ahler
assumption on the whole complex homogeneous manifold $G/\Gamma$ is investigated and characterized.\\

\vskip 2ex\noindent 
In contrast to this condition we study in the present paper the weaker "locally K\"ahler" property.
To be precise
 let $Y:=G_0/\Gamma \subset G/\Gamma =: X$ be the natural
inclusion of the Levi-flat compact homogeneous CR-solvmanifold $Y$
in the complex-homogeneous solvmanifold $X$.  We characterize
all the CR-solvmanifolds $Y$ which admit an embedding in some K\"ahler 
manifold. This is equivalent to assuming 
that there is an open subset $U \subset
G/\Gamma$ which is K\"ahler.
The main result that we prove is a characterization of this situation.  
It turns out that one of
the main necessary  conditions is that the restriction of the adjoint representation of 
$\g_0$ to $\m$ has imaginary spectrum and is diagonalizable.  
Note in passing that because of the set up, the manifold $Y$ is the base of 
the vector bundle in Mostow's conjecture.  

\vskip 2ex\noindent  
In addition to recalling some details from \cite{OR} that we need and proving the main result, 
we also present four examples to illustrate the theory, with two of these being locally K\"{a}hler 
and the other two not.   
The paper closes with a classification of locally K\"{a}hler compact CR-solvmanifolds of 
codimension one or two: a finite covering splits as a direct product of a Cousin group and some appropriate 
number of ${\mathbb C}^{*}$'s. 
This extends  a result in \cite{R} for the CR-hypersurface case, which is included due to the scope 
of the result and the method of its proof.

\section{The nilpotent case}  \label{nilcase}  

{\it For the rest of the paper we call
a triple $(\Gamma, G_0, G)$ consisting of a connected simply-connected solvable complex Lie group $G$, a real
connected (simply-connected) subgroup 
$G_0$ of $G$, and a discrete cocompact subgroup $\Gamma$ of $G_0$, a {\bf CR-solvmanifold}, or {\bf CRS} for short.  
We also assume throughout that $G_0/\Gamma$ is a generic CR-manifold of the complex manifold $G/\Gamma$, i.e. 
for the Lie algebras $\g_0$ and $\g$ of $G_0$ and $G$ one has
$\mathfrak{g}_0 + i \mathfrak{g}_0 =\mathfrak{g}$.}\\

If the CRS $(\Gamma, G_0, G)$ is {\it locally K\"ahler},
i.e., if $Y$ admits an open K\"ahler neighbourhood in $X$,
 which one may assume to be right-$G_0$-invariant with 
 right-$G_0$-invariant K\"ahler form $\omega$ (this is obviously equivalent
to saying that $Y$ is CR-embeddable in some K\"ahler manifold),    
then $\m$ is an abelian ideal in $\g$ and in $\g_0$
(see \cite{OR}, Corollary 7, p. 406).
Let $N$, $N_0$, $\mathfrak{n}$, $\mathfrak{n}_0$  
denote the associated nilradical objects. It is clear that $\m \subset \n$.   
A result of Mostow (\cite{Mos}) assures the existence of 
the so-called {\it nilradical fibration} which we recall in the next lemma.

\begin{Lem}\label{nilfibration}
The $N_0$-orbits (respectively the $N$-orbits) in $Y$ (respectively in $X$) are closed
and there is the following commutative diagram of fiber bundles
$$ 
    \begin{array}{ccc}  Y= G_0/\Gamma & \hookrightarrow & X= G/\Gamma  \\  
   \downarrow        &     & \downarrow \\  
     G_0 / N_0 \cdot \Gamma & \hookrightarrow & G / N \cdot \Gamma     \end{array}   
$$
the right vertical arrow being holomorphic.

\end{Lem}

\noindent 
It will be convenient if we now recall some results and remarks 
from section 3 in \cite{OR} that we need later.    
Since it is not always true that every subalgebra 
of a solvable Lie algebra has a complementary subalgebra, 
e.g., the center of a Heisenberg algebra has no complement,  
we first note the following in the locally K\"ahler CRS setting.    
In passing, we recall that Malcev \cite{Mal} showed that every Lie algebra 
over $\C$ admits a faithful representation into a splittable Lie algebra; 
see \cite{Re}  for another approach.    

\begin{Prop}\label{complement}  
Let $(\Gamma, G_0, G)$ be a locally K\"ahler CRS.    
Then there exists a complementary subalgebra $\a_0$ to $\m$ in $\g_{0}$, i.e., 
one has 
\[  
      \g_0 = \a_0  \ltimes  \m .  
\]  
\end{Prop}
{\tt Proof.} 
Let $$\mathfrak a_0=\{X \in \mathfrak g_0 \mid\omega(X,\m) = 0\}.$$
Then $\a_0$ is a real subalgebra of $\g_0$ and we have $\g_0 = \a_0  \ltimes  \m$.
 To see this 
consider the hermitian metric $h$ defined by $\omega$ and let $(v_1,...v_{2k},w_1,...,w_l)$ 
be an orthonormal basis of $\g_0$, 
where $(v_1,...v_{2k} )$ is an orthonormal basis of $\m$, such that 
$v_{j+k} = Jv_j$ for $1\leq j \leq k$, and $J$ denotes the complex structure.  

\vskip 2ex\noindent  
Using the relation $h(v,w)=\omega(Jv,w)$, it is easy to see that $\a_0$ is just the real 
span of $(w_1,..., w_l)$.   
Thus $\g_0= \a_0 \oplus  \m$ as a vector space. 
The fact that $\a_0$ is a Lie algebra follows from the formula 
$$ 0 = \omega(X, [Y, Z]) + \omega(Z, [X, Y]) + \omega(Y, [Z, X]), \leqno{(\diamond)}$$
for all $X, Y, Z \in  \g_0$, using the fact that $\m$ is an ideal in $\g_0$. Note that 
$[\g, \m] = [\g_0,\m]$. 
Furthermore $\m$ is abelian, so we 
get $[\g,\m] = [\a_0,\m]$. \qed

\begin{Th} [Theorem 2', \cite{OR}] \label{OR1988}
 Let $(\Gamma, G_0, G)$ be a locally K\"ahler CRS.   
 Let $\mathfrak z$ be
the center of $\mathfrak g$ and $\omega$ a right $G_0$-invariant K\"ahler form 
in a right $G_0$-invariant open neighbourhood of $G_0$ in $G$. 
Then $$ \mathfrak z \cap \m =\{X \in \m \mid \omega(X,[\mathfrak g, \m])=0\}.$$
In particular, if $\z = 0$, then $[\g,\m] = \m$.
\end{Th}
{\tt Proof.} 
 
By $(\diamond)$, it follows now  that 
\[   
        \z \cap \m \subset   \{X\in \m \mid \omega(X,[\g,\m]) = 0\} . 
\]  
To verify the opposite inclusion it is enough to show that for\\  $X_0 \in
   \{X \in \m \mid   \omega(X,[\g,\m]) = 0\}$  one has $[X_0,\a_0] = 0$. 
   For this let $Y \in \a_0$ and $Z \in \m$.  
By $(\diamond) $ we have 
\[   
      0=\omega(X_0, [Y, Z]) + \omega(Z, [X_0, Y]) + \omega(Y, [Z, X_0])  .
\]   
Since $[Y, Z]\in [\g,\m]$ and $[Z, X_0] \in \m$, it follows that 
\[  
         \omega(Z,[X_0,Y]) = 0 
\]  
for all $Z \in \m, \ Y  \in \a$.
But $X_0 \in \m$ implies that $[X_0,Y] \in \m$.  
Since $\m$ is a complex subspace of $\g$, 
we know that $\omega\vert_\m$ is nondegenerate. 
Therefore we have $[X_0,Y] = 0$  
Since the last remark is clear, the theorem is proved.  \qed

\begin{Cor}\label{Cor1}
Let $(\Gamma, G_0, G)$ be a locally K\"ahler CRS where $G$ is a nilpotent
complex Lie group.
Then $\m \subset  \z$, $\m \cap \g' = \{0\}$ and $\z = \z_0 + i\z_0$, where $\z_0$ 
is the center of $\g_0$.
\end{Cor}
{\tt Proof.} Firstly, it is clear that $[\m,\g]$ is an ideal in $\g$.
The preceding theorem gives an ideal splitting 
$$\m = (\m \cap \z) \oplus [\m,\g].$$
Consequently $[\g,[\m,\g]]=[\m,\g]$.  
For nilpotent $\g$ this is only possible if $[\m,\g]=0$, i.e. $\m \subset \z$.\\
This in turn implies that
$\g_0'=[\g_0,\g_0]=[\a_0,\a_0]\subset \a_0$. 
Therefore $\g_0' \cap \m =\{0\}$, $\g_0' \cap i \g_0' = 0$ and $\g' \cap \m=0$.\\
For the last assertion we first remark that
$\g_0$ being a generic subalgebra of $\g$, the center $\mathfrak z_0$ of 
$\mathfrak g_0$ is given by $\mathfrak z_0 = \mathfrak z \cap \mathfrak g_0$. 
Hence $\z_0 + i \z_0 \subset \z$.\\ 
Furthermore take $Z = X+iY \in \z,\ X,Y \in \g_0$. For an arbitrary
$X' \in \g_0$ we get
$$0=[Z,X']=[X,X']+i[Y,X'],$$
which shows that $[X,X'],[Y,X'] \in \g_0'\cap i\g_0'=0$. So $X,Y \in \z_0$
and $   \z \subset \z_0 + i \z_0$. \qed
 
\begin{Cor}\label{Cor1a}  
Let $(\Gamma, G_0, G)$ be a locally K\"ahler CRS where $G$ is a nilpotent
complex Lie group and suppose that there are only constant holomorphic
functions on the homogeneous manifold $X=G/\Gamma$. 
Then $G$ is an abelian complex Lie group.
\end{Cor}
 
{\tt Proof.} By hypothesis, $\mathcal O(G/\Gamma) = \mathbb C$.
By a result of \cite{BaOt}, we have an equivariant
holomorphic bundle 
\[  
      G/\Gamma \ra G/Z \Gamma ,  
\]
$Z$ being the center of $G$.   
Since the complex Lie group $M$ corresponding
to the Lie algebra $\m$ is contained in $Z$ by the previous corollary,
the homogeneous manifold $G/Z \Gamma$ is Stein, see \cite{GH1},
hence a point and $G=Z$ is abelian. \qed
 
\begin{Cor}\label{Cor2} 
Let $(\Gamma, G_0, G)$ be a CRS where $G$ is a nilpotent complex Lie group.   
Then the following three conditions are equivalent:\\
 1) $(\Gamma, G_0, G)$ is a locally K\"ahler CRS.\\
 2) $(\Gamma, G_0, G)$ is a K\"ahler CRS.\\
 3) There is a Lie algebra splitting $$\g_0=\a_0 \oplus \m,$$
where $\m$ is abelian and $\a_0 \cap i\a_0 =0$.
\end{Cor}
 
{\tt Proof.} 2)$\ \Rightarrow \ $1) is evident.\\
1)$\ \Rightarrow \ $3): In the proof of theorem \thesection.\ref{OR1988}
we have constructed the splitting and corollary \thesection.\ref{Cor1}
gives the additional properties.\\
3)$\ \Rightarrow \ $2): 
Assume that there is a splitting 
\[  
      \g_0=\a_0 \oplus \m  ,  
\]   
with $\m$ abelian and $\a_0 \cap i\a_0 =0$.  
Let $\a:=\a_0 \oplus  i\a_0$ and $A_0$, $A$ the corresponding Lie groups.  
Then we have
\[   
    G_0 = A_0 \times M, \ \ \  G = A \times M .   
\]  
Moreover there is a right $A_0$-invariant K\"ahler form on $A$, see \cite{Lo}. 
Since $M$ is abelian, this yields a right $G_0$-invariant K\"{a}hler form  
on $G$ and hence a K\"ahler form on $G/\Gamma$. \qed

%%%%%%%%%%%%%%%%%%%%%%%%  
 
\section{The Characterization}\label{GC}  
Before coming to the main theorem, we first need three lemmas.
 
\begin{Lem}\label{Regina}
Let $ \beta \in \mathbb C$ and $G:=\mathbb C \ltimes \mathbb C$
the complex Lie group with multiplication
\[   
    (z_0,w_0)\cdot (z,w):=(z_0+z,e^{\beta z_0} w + w_0).
\]    
Let further $G_0 := \mathbb R \ltimes \mathbb C \subset G$
be the subgroup of $G$ with $z \in \mathbb R$.    
Let $U$ be any connected left $G_0$-invariant open neighbourhood of $G_0$.   
Then there is a left $G_0$-invariant K\"ahler form
on $U$ if and only if $\Re e\ \beta =0$. The same conclusion holds if we replace
the left-invariant objects by right-invariant ones.
 \end{Lem}

{\tt Proof.} We explicitly write the left $G_0$-action on $U$:  
\[   
      (t,w_0)\cdot (z,w):=(t+z,e^{\beta t} w + w_0)   
\]  
 and let $H$ denote the complex Lie subgroup 
 $0 \ltimes \mathbb C \subset  G_0 \subset G$. \\
Let 
\[   
      \omega=f_1 dz \wedge d \bar{z}+f_2 dz \wedge d \bar{w}
      + f_3 d\bar{z} \wedge d w+f_4 dw \wedge d \bar{w}   
\]  
be a left $G_0$-invariant K\"ahler form on $U$.  
 
\vskip 2ex\noindent 
For $g_0=(t,w_0) \in G_0$ we get, using the
left $H$-invariance of $\omega$,  
\[   
    {g_0}^{\star}(\omega)=f_1(t+z ) dz \wedge d \bar{z}+
    e^{t \bar{\beta}}f_2(t+z ) dz \wedge d \bar{w}+$$ $$+
    e^{t \beta}f_3(t+z)  d\bar{z} \wedge 
    d w+ e^{t (\beta +\bar{\beta})}f_4(t+z ) dw \wedge d \bar{w} = \omega
\]   
for all $g_0 \in G_0$.  
But $d\omega=0$ now forces $f_4$ to be constant
which is the case if and only if $\Re e\ \beta =0$.\\
If $\Re e\ \beta =0$, one can take $\omega = dz \wedge d\bar{z}+
dw \wedge d\bar{w}$.\\
The composition with the inverse map in $G$ exchanges left and right-invariant
object and the last assertion is therefore evident.\qed  
\begin{Lem}\label{L1}
Let $\s_0$ be a real abelian Lie algebra ,
$\s = \s_0 \oplus i\s_0$ and $\t$ a complex abelian Lie algebra.   
Suppose that $\g_0= \s_0 \ltimes \t \subset \g= \s \ltimes \t$ is a semi-direct product 
such that the adjoint action of $\s_0$ on $\t$ has purely imaginary spectrum
and is diagonalisable. 
Then there is a right $G_0$-invariant K\"ahler form on $G$.
\end{Lem}
{\tt Proof. } By hypothesis,  the group
multiplication in the associated simply-connected
complex Lie group $G$ is given in suitable  coordinates by
 \[  
             ( (s_1,...,s_n) , (t_1,...,t_m) ) \circ  ( (s^0_1,...,s^0_n) , (t^0_1,...,t^0_m) ) =
 \]   
 \[   
              = ((s_1+s^0_1,...,s_n+s_n^0),(e^{i\sum_{j=1}^n\alpha_{j1}  s_j} t_1^0+ 
                     t_1,..., e^{i\sum_{j=1}^n\alpha_{jm}  s_j} t_m^0+ t_m)) , 
 \]   
with real constants $\alpha_{jk}$.
The K\"ahler form $\omega := \sum ds^0_j \wedge d{\bar s^0}_j
+ \sum dt^0_k \wedge d{\bar t^0}_k$ on $G$ is then obviously left $G_0$-invariant. 
By composing with the inverse map, one produces in a standard way a
right $G_0$-invariant K\"ahler form on $G$.
 \qed
 
\begin{Lem}\label{L2}
Suppose $(\Gamma,G_0, G)$ is a  
totally real CRS, i.e. $\g_0 \cap i\g_0 =\m = \{0\}$.
Then $(\Gamma, G_0, G)$ is a locally K\"ahler CRS. 
Furthermore, for any connected complex Lie subgroup 
$H \subset G$ set $H_0 := G_0 \cap H$.  
Then there is a right $H_0$-invariant neighborhood 
$U$ of $H_0$ in $G$ and a right $H_0$-invariant K\"ahler form on $U$.    
\end{Lem}

{\tt Proof. } 
Since $G_{0}/\Gamma$ is a totally real, real analytic submanifold 
of $G/\Gamma$, there exists an open neighbourhood $U$ of $G_{0}/\Gamma$ 
in $G/\Gamma$ that is Stein \cite{Gr}.  
By shrinking and renaming, if necessary, we may assume that $U$ is 
right $G_{0}$-invariant.  
It is clear that a K\"{a}hler form exists on $U$.    
Since $G_0/\Gamma$ is compact, we may average the K\"{a}hler form, as in the proof of 
Proposition 1 (see p. 166) in \cite{GOR}, and so obtain a form that is right $G_{0}$-invariant.    
The remaining statements follow from this.  \qed

\begin{Def}
We call a pair $(H_0,H)$ as in the preceding lemma a {\it locally K\"ahler pair}.
\end{Def}
 
In this section we shall prove the following
 
\begin{Th}\label{MT}
Let $(\Gamma, G_0, G)$ be a CRS. Then the following  
statements are equvalent:\\
1) $(\Gamma, G_0, G)$ is a locally K\"ahler CRS.\\
2) $(\Gamma, G_0, G)$ satisfies the conditions
\begin{itemize}
\item[]{i)}  $(\Gamma \cap N_0, N_0, N)$ is a locally K\"ahler CRS
\item[]{ii)} there is a Lie algebra splitting $\g_0 =\a_0 \ltimes \m$
\item[]{iii)}  the adjoint representation $ad(\g_0)\vert_{\m}$ of $\g_0$ on
$\m$ is diagonalisable and has purely imaginary spectrum.
\end{itemize}
\end{Th}
{\tt Proof. }  1)$\, \Rightarrow\, $2):\\ 
Suppose that $(\Gamma, G_0, G)$ is a locally K\"ahler CRS.
By restriction it follows immediately that 
$(\Gamma \cap N_0, N_0, N)$ is a locally K\"ahler CRS, i.e. 2) i).    
Condition 2) ii) is also
necessary as we have seen in proposition 4.\ref{complement}.
By lemma 4.\ref{Regina} we get the purely imaginary spectrum of 
$ad(\g_0)\vert_{\m}$ and by the methods of \cite{GOR},
p.p. 166-167 one can prove the diagonalizability in condition 2) iii), as we 
now point out.    
Since the K\"{a}hler form $\omega$ restricted to $\g_0$ may be assumed 
to be right invariant and $d\omega =0$, one has 
\[ 
      \omega(X,[Y,Z]) \; + \; \omega(Z,[X,Y]) \; + \; \omega(Y,[Z,X]) \; = \; 0 
\]   
for all $X,Y,Z \in \g_0$.  
Now suppose $v$ is an eigenvector for the adjoint action on $\m$ and 
set $v^{\perp} := \{ w\in\m \ | \ \omega(w, \langle v \rangle ) = 0 \}$. 
Then for $w \in v^{\perp}$ and $X \in \g_0$ one has  
\[  
      \omega(X,[w,v]) \; + \; \omega(v,[X,w]) \; + \; \omega(w,[v,X]) \; = \; 0  .  
\]   
Now $\m$ is abelian, so $[w,v]=0$.  
Also $[v,X]= \mu v$ for some $\mu\in\C$, so 
$\omega(w,[v,X])=\mu\omega(w,v)=0$.  
It follows that $[X,w]\in v^{\perp}$.  
Hence, by induction, the adjoint action of $\g_0$ restricted to $\m$ is diagonalizable.   
 
\vskip 2ex\noindent  
2)$\, \Rightarrow\, $1):\\ 
Now assume the three conditions in 2).   
We define $\h_0:=\a_0  \cap \n$ and $\h:=\a  \cap \n$. 
Then $\h= \h_0 \oplus i \h_0$ as a vector space.   
Furthermore $\n_0 =\h_0 \oplus \m$ and 
$\n =\h \oplus \m$ as Lie algebras in view of condition 2) i) and 
Corollary 2.4, which also gives that $\h_0\,  \triangleleft\,  \g_0$ and 
$\h\,  \triangleleft\,  \g$.  
Let $H_0 \tl G_0$ and $H \tl G$ denote the subgroups of
$G$ corresponding to $\h_0$ and $\h$ and $\pi_1:G\ra G/H$ the projection. 
We then have that  
\[   
        L_0:=G_0/H_0 \simeq  A_0/H_0 \ltimes M=: B_0 \ltimes M    
\]   
and 
\[    
        L := G/H \simeq A/H \ltimes M =: B \ltimes M   
\]   
are Lie groups.   
The new Lie algebras are denoted by
the corresponding frakture letters.   
One has $\l_0 \cap i\l_0 =\m,\ \l_0 + i\l_0 = \l$ and condition 2) iii) now  
gives that ${\rm ad}(\b_0)\vert_{\m}$ has purely imaginary spectrum 
and is diagonalizable. 
So by lemma \thesection.\ref{L1} there is a 
right $L_0$-invariant K\"ahler form $\omega_1$ on $L$ whose pullback to 
$G$ will also be denoted by $\omega_1$ and which is right $G_0$-invariant. 

\vskip 2ex\noindent   
Define now ${\tilde \g}:= \g_0 \oplus i\g_0 $ the ``formal'' complexification of
$\g_0$ and consider the associated totally real CRS $(\Gamma, G_0, \tilde G)$
which is locally K\"ahler by lemma \thesection.\ref{L2}.    
This lemma shows furthermore that the pair $(A_0,A)$ is locally K\"ahler.   
Let $\omega_2$ denote the local right $A_0$-invariant K\"ahler
form on $A$  and remark that dividing the pair $(G_0,G)$ by the complex
Lie group $M$ gives exactly the pair $(A_0,A)$. 
Pulling back $\omega_2$ to $G$ we again get a right $G_0$-invariant form.   

\vskip 2ex\noindent  
Finally the sum $\omega:= \omega_1+ \omega_2$ is non-degenerate and therefore is a right 
$G_0$-invariant K\"ahler form on a right $G_0$-invariant open neighborhood of $G_0$ in $G$.
This proves that $(\Gamma, G_0, G)$ is a locally K\"ahler CRS. \qed

%%%%%%%%%%%%%%%%%%%%%

\section{Examples}  \label{Exs}  

We now present four examples of CRS's $(\Gamma,G_{0},G)$
that illustrate various aspects of the general theory.  
The second and third of these are locally K\"{a}hler, while the other two are not.  

\subsection{Non-imaginary Spectrum}  

This example fibers as a 2-dimensional Cousin group bundle over $\C^*$.   
The critical point is the existence of an eigenvalue that is not purely imaginary.  
This example is modeled on a construction of Inoue surfaces, for example see \cite{BPV}.  

\vskip 2ex\noindent 
We choose a matrix $A \in SL_{3}(\Z)$ that has one real eigenvalue $\alpha>1$ 
and two nonreal complex conjugate ones $\beta$ and $\overline{\beta}$.  
To be precise one may pick 
\[ 
        A \; = \; \left( \begin{array}{ccc}  0 & 1 & 0 \\ k & 0 & 1 \\ 1 & 1-k & 0 \end{array} \right) 
\]   
where $k$ is an integer.  
Further, we let $(a_{1},a_{2},a_{3})$ be a real eigenvector corresponding to $\alpha$ 
and $(b_{1},b_{2},b_{3})$ an eigenvector corresponding to $\beta$.

\vskip 2ex\noindent 
Consider the group $G := \C \ltimes \C^2$ with group structure given by 
\[  
     (t,x,y) \cdot (t_{0},x_{0},y_{0}) \; := \; (t+t_{0}, x + e^{t\log\alpha}x_{0}, y+e^{t\log\beta}y_{0}),
\]   
where $\log $ is the principal branch of the logarithm on $\mathbb C \setminus \mathbb R^+$.
 The discrete subgroup $\Gamma$ is generated by the following elements:  
\[  
     \Gamma \; := \; \langle (t, 0 , 0), (0, a_{1},b_{1}), (0, a_{2},b_{2}),(0, a_{3},b_{3}) \rangle_{\Z}  
\]  
where  $t\in\Z$.  
Clearly, $G' = \{ (0, x, y) | x,y\in\C\} $, has closed orbits in $G/\Gamma$,  
and the fibration $G/\Gamma \to G/G'\cdot\Gamma$ is a 
Cousin group bundle over $\C^*$.  
Note that for any choice of the branch of the logarithm,  
$\log\beta$ is purely imaginary if and only if $|\beta|=1$.  
However, this is not possible, because $\alpha\beta\overline{\beta}=\det A = 1$ and $\alpha > 1$.  
The homogeneous space $G/\Gamma$ is not K\"{a}hler, and 
Lemma 4.\ref{Regina} shows that this example is not locally K\"ahler either.  
 
%%%%%%%%%%%%%%%%%%%%%%%    

\subsection{A $5$-dimensional locally K\"ahler example}  
Let $H_{3}$ be the 3-dimensional Heisenberg group and let the complex Lie group 
${\rm GL}(2,\C)$ act as a group of holomorphic  transformations on $H_{3}$,     
where  for $\displaystyle  A \; = \; \left( \begin{array}{cc} a & b \\ c & d \end{array} \right)   \in {\rm GL}(2,\C)$
the action is given by    
\[ 
        \left( \begin{array}{ccc}  1 & x & z \\ 0 & 1 & y \\ 0 & 0 & 1 \end{array} \right) 
        \; \stackrel{F_{A}}{\mapsto} \; 
 \left( \begin{array}{ccc}  1 & ax+by & (\det A)\left[ z -xy/2 \right] + 
 (ax+by)(cx+dy)/2  \\ 
 0 & 1 &cx+dy \\ 0 & 0 & 1 \end{array} \right) 
\]  
A direct calculation shows that the map $A \mapsto F_{A}$ is a group homomorphism.   

\vskip 2ex\noindent 
The subgroup 
\[ 
    \Lambda \; := \; \left\{ \left( \begin{array}{ccc}  1 & n & k/2 \\ 0 & 1 & m \\ 0 & 0 & 1\end{array} \right) 
    \left|\;  n, m, k \in \Z    \right.  \right\} 
\]  
is discrete in $H_{3}$ and note that for $\displaystyle 
    A \; = \; \left( \begin{array}{cc} 2 & 1 \\ 1 & 1 \end{array} \right) $  
one has $F_{A}(\Lambda) = \Lambda$.  
We now let 
\[ 
B \; = \; \left( \begin{array}{cc} 1 & (\sqrt{5}-1)/2 \\ -1 & (\sqrt{5}+1)/2 \end{array} \right) 
\]  
Then $F_{B}(\Lambda)$ is a discrete subgroup of the Heisenberg group.    

\vskip 2ex\noindent  
Next consider the representation from $\C$ into the group of automorphisms of $H_{3}$, where 
$t\in\C$ maps to conjugation by the matrix 
\[ 
      \left( \begin{array}{ccc} 1 & 0 & 0 \\ 0 & e^{t} & 0 \\ 0 & 0 & 1 \end{array}\right) 
\]  
Explicitly, this conjugation is as follows 
 \[ 
      \left( \begin{array}{ccc} 1 & 0 & 0 \\ 0 & e^{-t} & 0 \\ 0 & 0 & 1 \end{array}\right) 
      \left( \begin{array}{ccc} 1 & x & z \\ 0 & 1 & y \\ 0 & 0 & 1 \end{array}\right) 
      \left( \begin{array}{ccc} 1 & 0 & 0 \\ 0 & e^{t} & 0 \\ 0 & 0 & 1 \end{array}\right) 
      \; = \;  \left( \begin{array}{ccc} 1 & e^{t}x & z \\ 0 & 1 & e^{-t}y \\ 0 & 0 & 1 \end{array}\right) 
 \]  

\vskip 2ex\noindent  
Now let $K=\R$ or $K=\C$ and $(t,x,y,z) \in K^4$. 
Define a group structure $G_K$ on $K^4$ by  
 \[   
     (t,x,y,z)\circ (t',x',y',z'):= (t+t',x+e^t x', y+e^{-t} y', z+z'+ xy'e^{-t}).  
 \]     
 
 \vskip 2ex\noindent 
Let $\lambda:= \frac{1}{2} (3 + \sqrt{5})$ and $\alpha:= \ln \lambda$.
Let $\Gamma$ be the subgroup of $G_{\R}$ generated by the elements 
\[  
   (\alpha,0,0,0); (0,1,-1,- \frac{1}{2}); (0,  \frac{1}{2}(\sqrt{5}-1),\frac{1}{2}(\sqrt{5}+1), \frac{1}{2}); (0, 0, 0, \sqrt{5}) .  
\]  
The remarks above imply that $\Gamma$ is discrete.  
Further $\Gamma$ is a cocompact subgroup of $G_{\R}$.  

\vskip 2ex\noindent 
Finally, define the 5-dimensional solvable complex Lie group 
$G := G_{\C}\times \C$, where $G_{\C}$ is embedded as the first component. 
We add to the embedded $\Gamma$ the following two generators 
\[ 
      (0,0,0,0,1) \quad \mbox{ and }\quad (0,0,0,\sqrt{2}, i \sqrt{3})  
\]  
and this gives a Cousin group structure in the last two coordinates.  
The subgroup $G_{0}$ is defined to be the real span of the discrete subgroup.    

\vskip 2ex\noindent 
One should note that in this example the complex ideal $\mathfrak m$ is contained 
in the product of the center of the Heisenberg group and the trivial $\C$-factor.  
Therefore the adjoint action of $\mathfrak g_{0}$ restricted to $\mathfrak m$ is trivial.    
As a consequence of theorem 5.\ref{MT} we see that $(\Gamma, G_0, G)$ is
a locally K\"ahler CRS.

%%%%%%%%%%%%%%%%%%%%%%   

\subsection{A $4$-dimensional locally K\"ahler example}  \label{Salem}  

In the previous example the adjoint action restricted to $\mathfrak m$ is trivial.  
We would next like to present an example where this action is purely imaginary, 
but in a nontrivial way, i.e., the action has a nonzero eigenvalue.      

\vskip 2ex \noindent  
Let $p(x):=x^4-3x^3+3x^2-3x+1$ 
(one can also take $x^4-x^3-x^2-x+1$ , by the way).
Then $p$ is $\mathbb Q$-irreducible and the four roots
of $p$ are $\alpha >1$, $0< \alpha^{-1}< 1$, $\alpha \in \mathbb R$, 
$\beta = e^{is}$ and $\bar{\beta}$ with $\vert \beta \vert =1$, i.e.
$\beta$ is on the unit circle and is multiplicatively of infinite order, 
i.e. the group $\{ {\beta}^k \mid k \in \mathbb Z \} $ is isomorphic to
$\mathbb Z$.  In fact, $\alpha$ is a unit in the ring of algebraic integers.
The number $\alpha$ is, by definition, a {\it Salem number}.   

\vskip 2ex\noindent 
Now let $K:=\mathbb Q[\alpha]$ be the number field of degree
$4$ and $O_K \simeq \mathbb Z^4$ the ring of integers of $K$.
Let $\sigma_1, \sigma_2 : K \rightarrow \mathbb R$
be the real imbeddings  and $\sigma_3, \sigma_4=\bar{\sigma_3}:
K \rightarrow \mathbb C$ the complex, non-real imbeddings of $K$.
Define $\sigma : K \rightarrow \mathbb C^3$ by
\[   
          \sigma(k):= (\sigma_1(k),\sigma_2(k),\sigma_3(k)) .  
\]  
It is easy to check that $\Lambda:= \sigma(O_K)\simeq \mathbb Z^4$ 
is discrete in $ \mathbb C^3$ and that the quotient $\C^3/\Lambda$ is a Cousin group; 
see \cite{OT}.
Let also   
\[  
 D:=\left(\begin{array}{ccc}\ln \alpha & 0 & 0 \\ 0 &-\ln \alpha & 0  \\
0 & 0 & is \end{array}\right).  
\]   
Define now a solvable $G$ group structure on $\mathbb C \ltimes \mathbb C^3$
by
\[   
        (z,b)(z_0,b_0) := (z+z_0, e^{zD}(b_0) +b)  .  
\]  
Let $V_{\mathbb R}$ be the real span of $\Lambda$ in $\mathbb C^3$
and $G_0:= \mathbb R \ltimes V_{\mathbb R}$ be the corresponding real subgroup of $G$. 
Then $\Gamma:=\mathbb Z \ltimes \Lambda $
is a discrete cocompact subgroup of $G_0$ and $(\Gamma, G_0, G)$
is a CRS.   

\vskip 2ex\noindent  
Further note that the subgroup $M:=V_{\mathbb R} \cap i V_{\mathbb R}$
has dense orbits in $V_{\mathbb R}/\Lambda \simeq (S^1)^4$ and that the 
Lie algebra $\mathfrak g_0$ has purely imaginary spectrum on the Lie
algebra $\mathfrak m$ of the group $M$.  
In fact, in our realisation $M= 0 \times 0 \times \mathbb C \subset \mathbb C^3$!  
Since $\mathfrak g_0$ does NOT have a purely imaginary spectrum on itself,   
$(\Gamma, G_0, G)$ is NOT a K\"ahler CRS.
But the following holds.  

\begin{Lem}
 $(\Gamma, G_0, G)$ is a  locally K\"ahler CRS.
\end{Lem}

{\it Proof.}   All conditions in 2) of theorem 5.\ref{MT} are satisfied,
therefore $(\Gamma, G_0, G)$ is a locally K\"ahler CRS. \qed 
%%%%%%%%%%%%%%%%%%%%%%%%  
%%%%%%%%%%%%%%%%%%%%%%%%  

\subsection{Non-diagonalizable Example}

We now show by an example that the diagonalizability assumption is needed.  

\vskip 2ex\noindent 
Define the real $7$-dimensional solvable Lie group $G_0$ as the semi-direct
product of $(\mathbb R,+) $ and $(\mathbb C^3,+)$ with multiplication
\[  
    (t,z, \left(\begin{array}{c} a \\  b\end{array}\right))\circ 
   (t_0,z_0,  \left(\begin{array}{c} a_0 \\  b_0\end{array}\right))
   := (t+t_0, e^{2\pi i t} z_0 + z, e^{tD}  \left(\begin{array}{c} a_0 \\  b_0\end{array}\right) +  
   \left(\begin{array}{c} a \\  b\end{array}\right))    ,
\]   
where $D:=  \left(\begin{array}{cc}0  & 1  \\  0& 0 \end{array}\right)$.

\vskip 2ex\noindent 
Let $G$ be the complex Lie group given by the corresponding semi-direct product
of $(\mathbb C,+) $ and $(\mathbb C^3,+)$. 
Now let $\Gamma$ be the discrete subgroup of $G_0$ given by the semi-direct product of
$(\mathbb Z,+) $ and $(  (\Z + i\Z)^3,+)$. 
Then it is easy to see that $\Gamma$
is cocompact in $G_0$. In this example property 2) of Theorem 5.\ref{MT} 
is satisfied except for diagonalizability.   
Furthermore example 6b), p.414 in \cite{OR} shows that $Y$ is not locally K\"{a}hler in $X$.

%%%%%%%%%%%%%%%%%%%   

  %%%%%%%%%%%%%%%%%%%%%%%%%%%%%%%% 
\section{A Remark on Holomorphic Reductions}  

\noindent  
We make some simple observations in this short section
about holomorphic reductions in the locally K\"{a}hler setting.   
In \cite{OR} it is shown that the fiber of the holomorphic reduction of a K\"ahler CRS is a Cousin group.  
This is no longer the case, in general, for the globalization of  locally K\"{a}hler CR-solvmanifolds, 
but we can extend Corollary 3.\ref{Cor1a} in the following way.

\begin{Prop}  \label{noholofcns} 
Let $(\Gamma,G_{0},G)$ be a locally K\"ahler CRS.       
Then ${\cal O}(G/\Gamma) = {\mathbb C}$ if and only if $G/\Gamma$ is a Cousin group.     
\end{Prop}  

{\tt Proof.} 
By lemma \ref{nilcase}.\ref{nilfibration}, we may consider the nilradical fibration  
$G/\Gamma \to G/N\cdot\Gamma$ of the complex solvmanifold $X = G/\Gamma$.  
Its base $G/N\cdot\Gamma$ is a Stein abelian Lie group, 
because $G' \subset N$ and ${\mathfrak m} \; \subset \; {\mathfrak n}$.  
Thus $G=N$ is nilpotent whenever ${\cal O}(G/\Gamma) = {\mathbb C}$.  
Since locally K\"{a}hler and K\"{a}hler are equivalent in the nilpotent setting, 
see Corollary \ref{nilcase}.\ref{Cor2}, one direction   
follows from Corollary \ref{nilcase}.\ref{Cor1a}.  
The other direction is evident.  
\qed 

\vskip 2ex\noindent  {\bf Remark:}\ 
There are always a finite number of fibrations given by successive holomorphic reductions.  
However, we would now like to note what happens in the locally K\"{a}hler setting.  
Suppose $G_0/H_0$ is a locally K\"{a}hler CRS with globalization $G/H$.   
Let 
\begin{eqnarray}\label{1stred}
         G/H \; \longrightarrow \; G/J_1
 \end{eqnarray}  
be its holomorphic reduction.    
Denote by $G_{N}(H^{\circ})$ the normalizer in $G$ of the connected component of the identity $H^{\circ}$ of $H$.  
The base $G/N_{G}(H^{\circ})$ of the normalizer fibration of $G/H$ is biholomorphic to $({\mathbb C}^{*})^{k}$ for some $k$, see \cite{GH2}.  
This implies that $J_{1}$ is a subgroup of $G_{N}(H^{\circ})$, the complex ideal $\m := {\mathfrak g}_0 \cap i  {\mathfrak g}_0$ is contained 
in the Lie algebra of $J_{1}$, and $J_{1}/H$ is parallelizable.  
Now set $\widehat{J_{1}} : = J/H^{\circ}$ and $\Gamma := H/H^{\circ}$ and let  $\widehat{\m} := \m/{\mathfrak h}\cap\m $ be the quotient of $\m$ 
by the Lie algebra $\mathfrak h$ of $H^{\circ}$ intersected with $\m$.  
Let  
\begin{eqnarray}\label{2ndred}   
          \widehat{J_{1}}/\Gamma   \; \longrightarrow \;   \widehat{J_{1}}/\widehat{J_{2}}  
\end{eqnarray}  
be the holomorphic reduction of the fiber of the fibration given in (\ref{1stred}).  
As observed in the proof of the previous Proposition, the base  $\widehat{J_{1}}/\Gamma\cdot N_{\widehat{J_{1}}}$ 
of the nilradical fibration of $\widehat{J_{1}}/\Gamma$ is a Stein abelian Lie group isomorphic to $({\mathbb C}^{*})^{l}$ for some $l$.  
This implies that $\widehat{J_{2}}$ is a subgroup of $\Gamma\cdot N_{\widehat{J_{1}}}$ and that the complex ideal $\widehat{\m}$ 
is contained in the Lie algebra of the group  $\widehat{J_{2}}$.  
Since the base $\widehat{J_{1}}/\widehat{J_{2}}$ is a holomorphically separable solvmanifold, it is Stein, see \cite{HO}.   
Thus the fiber $ \widehat{J_{2}}/\Gamma$ of the fibration in (\ref{2ndred}) is connected by a theorem of K. Stein.  
Now the nilpotent group $\widehat{J_{2}}^{\circ}$ acts transitively on the connected fiber of (\ref{2ndred}).    
By Corollary \ref{nilcase}.\ref{Cor2} the space $\widehat{J_{2}}/\Gamma$  is K\"{a}hler.  
Let $\widehat{J_{2}}/\Gamma \to   \widehat{J_{2}}/\widehat{J_{3}}$ be its holomorphic reduction.  
Its fiber $\widehat{J_{3}}/\Gamma$ is a Cousin group by \cite{OR}.    
Finally we have the following tower of fibrations  
\begin{eqnarray}\label{genholored}   
     G/H \; \stackrel{J_{3}/H}{\longrightarrow} \;  G/J_{3}  \; \longrightarrow \; G/J_{2} 
      \; \longrightarrow \; G/J_{1} ,  
\end{eqnarray}   
because of the fact that holomorphic reductions are independent of the particular group that is acting transitively 
on the space.

\vskip 2ex\noindent{\bf Caveat:}\  The space $G/J_3$ is not necessarily holomorphically separable, 
i.e., the fibration $G/H \to G/J_3$ need not be the holomorphic reduction of $G/H$. 
This is illustrated by the Coeur\'{e}-Loeb example, \cite{CL}, where the holomorphic reduction 
has $\C^*\times\C^*$ as fiber; i.e., one has the following 
\[  
     G/H \; = \;  G/J_{3}  \; \stackrel{{\C}^*\times{\C}^*}{\longrightarrow} \; G/J_{2}  \; = \; G/J_{1} \; = \; {\C}^* .    
\]   

%%%%%%%%%%%%%%%%%%%%%%%%%%%    

\vskip 2ex\noindent{\bf Remark:}\   
In Example \ref{Salem} involving the Salem numbers, one should note that  
because of the nature of the adjoint action of the Lie algebra 
${\mathfrak g}_{0}$ on $\m$, it follows that neither the first bundle given in (\ref{genholored}), 
namely $G/H \to G/J_{3}$,  
nor any finite covering of it can have the structure of a {\bf principal} Cousin bundle. 
Thus the locally K\"{a}hler case is in sharp contrast with the K\"{a}hler case;  
in \cite{GO} we showed that every K\"{a}hler solvmanifold has a finite covering 
whose holomorphic reduction is a principal Cousin bundle.

%%%%%%%%%%%%%%%%%%%%%%%%

\section{Classification in Low Codimensions}  

Some of the first results concerning compact, locally K\"{a}hler, 
homogeneous hypersurfaces appear in the dissertation of W. Richthofer \cite{R}.  
Using the tools at hand we now present a classification  
of locally K\"{a}hler CRS  $(\Gamma,G_{0},G)$ in codimensions one and two.  
  
\begin{Th} 
Suppose $(\Gamma, G_{0},G)$ is a locally K\"{a}hler CRS of codimension at most two.  
\begin{enumerate}
\item If $G_{0}/\Gamma$ is a hypersurface, then one of the following occurs:  
\subitem (i) $G/\Gamma$ is a Cousin group, or, 
\subitem (ii) A finite covering of $G/\Gamma$ is a product of a torus and $\C^*$.  
\item If $G_{0}/\Gamma$ has codimension two, then one of the following occurs:  
\subitem (i)  $G/\Gamma$ is a Cousin group, or, 
\subitem (ii) A finite covering of $G/\Gamma$ is a product of a Cousin group (of 
hypersurface type) and $\C^*$, or, 
\subitem (iii)  A finite covering of $G/\Gamma$ is a product of a torus and $\C^*\times\C^*$.  
\end{enumerate}  
\end{Th}  

 {\tt Proof.}  {\bf Codimension 1:}\  
Proposition 4.\ref{noholofcns} handles the case when ${\cal O}(G/\Gamma) = {\mathbb C}$.  
Suppose ${\cal O}(G/\Gamma) \not= {\mathbb C}$ and let 
\[  
     G/\Gamma \; \longrightarrow \; G/J\cdot\Gamma 
\]  
be the holomorphic reduction of $G/\Gamma$.  
For codimension reasons its base is biholomorphic 
to $\C^{*}$ and its fiber is a compact complex torus.  
It then follows from Proposition 1 in [GOR] that a finite 
covering of $G/\Gamma$ is biholomorphic to a product.  

\vskip 2ex\noindent{\bf Codimension 2:}\  
The case of no holomorphic functions is again handled by  proposition 4.\ref{noholofcns}.     
Using the same notation we again consider the holomorphic reduction.  
Its base cannot be $\C$, since $G_{0}/\Gamma$ would then be complex,  
contrary to the generic assumption.  
So its base is either $\C^{*}$ or $\C^{*}\times\C^{*}$.  
In the second case its fiber is compact and a finite covering splits as a product.  

\vskip 2ex\noindent 
One is reduced to considering the case where the base of the holomorphic reduction 
 \[  
    X= G/\Gamma \; \longrightarrow \; G/J\cdot\Gamma = \C^{*}
\]  
is $\C^{*}$ and its fiber $F := J/J\cdot\Gamma = J /(J \cap \Gamma)$ is a codimension 
one locally K\"{a}hler CRS.  
First we claim that $F$ is a hypersurface Cousin group.  
By the codimension one case considered above either $F$ is a Cousin group or 
there is a splitting $J / (J\cap\Gamma) = J/H \times H/(H\cap\Gamma)$, where 
we have abused the language and assumed that the space itself splits.  
Assume we are in the second case.  
If $H/(H\cap\Gamma)$ were biholomorphic to a torus, one would then 
have an intermediate fibration 
\[  
     G/\Gamma \; \longrightarrow \;  G/H \cdot\Gamma 
\]  
with fiber the torus and base a $\C^*$-bundle over $\C^*$.  
By the classification of two dimensional solvmanifolds a finite 
covering of $G/H\cdot\Gamma$ would be biholomorphic to 
a direct product $\C^*\times\C^*$ and this would contradict the 
fact that the base of the holomorphic reduction of $X$ has dimension one.  
Thus $F$ is a Cousin group.          
Now we claim that a finite cover of $X$ is the trivial bundle.      

\vskip 2ex\noindent  
Clearly, $J \subset N \subset G$ and $J$ has complex codimension one in $G$.  
So if $N=G$, then we are in the nilpotent setting and this case is easily handled as follows.  
The complex ideal $\mathfrak m$ lies in the center $\mathfrak z$ of the Lie 
algebra $\mathfrak g$ of $G$, see Theorem 2' in \cite{OR}.    
If $\mathfrak m = \mathfrak z$, then $J$ is the center of $G$ and has closed orbits 
and the base of the holomorphic reduction of $G/\Gamma$ has dimension two.  
(This is a consequence of the construction of the holomorphic reduction 
of complex nilmanifolds, see \cite{GH1}.)  
Therefore, $\mathfrak m$ has codimension at least one in $\mathfrak z$.  
But then the center of $G$ has codimension at most one and thus $G$ is abelian.  
The rest is obvious.  
 
 \vskip 2ex\noindent 
So we assume that $J=N$ and the fiber is a Cousin group of CR-hypersurface type. 
It is clear that $\rm{Ad}(\Gamma)  \subset \rm{GL}(J_0 \cap \Gamma)$,   
is diagonalisable (as an element of $\rm{GL}(J)$ over $\mathbb C$), 
and its restriction to $M$ is contained in the compact torus acting on $M$.  
Since it stabilizes a discrete group and is invertible, it has a generator of determinant $1$.   
Since $M$ is of codimension one in $J_0$, there is one   
remaining eigenvalue for the adjoint action of the Lie algebra of $J_{0}$ 
on itself and this also has absolute value $1$.   
But all this together implies that $\rm{Ad}(\Gamma)$ is a finite group.   
The rest of the proof is identical to the proof of proposition 1 in [GOR], p. 167.  
\qed

\vskip 2ex\noindent{\bf Remark:}\ 
One should note that this proof relies on the codimension two assumption.  
The example constructed using the Salem number (see \ref{Salem}) 
and the example on p. 413 
in \cite{OR} are both of codimension three and in both cases 
no finite covering of their holomorphic reductions splits as a product.  

%%%%%%%%%%%%%%%%%%%%%%%%    

\end{document}